\documentclass[11pt,a4paper]{amsart}
\usepackage{hyperref}
\usepackage{graphicx} 
\usepackage{wasysym}
\usepackage[mathscr]{eucal}

\usepackage{amsmath}

\usepackage{amssymb}
\usepackage{amsthm}
\usepackage{mathtools} 

\usepackage{amsfonts}
\usepackage{tikz}
\usepackage{enumerate}
\usepackage{csquotes}
\usepackage[draft]{changes}
\newtheorem{theorem}{Theorem}[section]

\newtheorem{proposition}[theorem]{Proposition}
\newtheorem{corollary}[theorem]{Corollary}

\theoremstyle{definition}
\newtheorem{definition}[theorem]{Definition}
\newtheorem{example}[theorem]{Example}

\newtheorem{remark}[theorem]{remark}

\providecommand{\ceil}[1]{{}^\lceil #1 {}^\rceil}

\begin{document}
	
	\title[New results about aggregation functions of quasi-pseudometric modulars]{New results about aggregation functions of quasi-pseudometric modulars}
	
	\author[A. Fructuoso-Bonet]{Alejandro Fructuoso-Bonet}
	\address[A. Fructuoso-Bonet]{Instituto Universitario de Matem\'atica Pura y Aplicada, Universitat Polit\`ecnica de Val\`encia, Camino de Vera s/n, 46022 Valencia, Spain}
	\email{afrubon@posgrado.upv.es}
	\author[J. Rodr\'{\i}guez-L\'opez]{Jes\'us Rodr\'{\i}guez-L\'opez}
	\address[J. Rodr\'{\i}guez-L\'opez]{Instituto Universitario de Matem\'atica Pura y Aplicada, Universitat Polit\`ecnica de Val\`encia, Camino de Vera s/n, 46022 Valencia, Spain}
	\thanks{The last author's research
		is part of the project PID2022-139248NB-I00 funded by MICIU/AEI/10.13039/501100011033 and ERDF/EU}
	\email{jrlopez@mat.upv.es}
	
	

	\begin{abstract}
		In recent studies, Bibiloni-Femenias, Mi\~{n}ana and Valero characterized the functions that aggregate a family of (quasi-)(pseudo)metric modulars defined over a fixed set $X$ into a single one. In this paper, we adopt a related but different approach to examine those functions that allow us to define a (quasi-)(pseudo)metric modular in the Cartesian product of (quasi-)(pseudo)metric modular spaces. We base our research on the recent development of a general theory of aggregation functions between quantales. This enables to shed light between the two different ways of aggregation (quasi-)(pseudo)metric modulars.  
	\end{abstract}
	
	\date{This is the author's version of the article accepted for publication in \textit{Mathematics}. 
		The final published version is available at MDPI via DOI: \href{https://doi.org/10.3390/math13050809}{10.3390/math13050809}.}
	\maketitle

%
\section{Introduction}
%

Aggregation refers to the process of gathering individual items and combining them into a unique one.
This action appears in mathematics for different types of objects. The most common are numbers, and the thought behind aggregating them is to obtain a representative one that summarizes the information in the set of numbers. Thus, in a general context, an aggregation function $F$ is of the form $F:X^n\to X$, where $X$ is a nonempty set. This simple idea has gained significant importance in a wide range of disciplines due to its ability to model decision-making problems (see, for example, \cite{Pap15,BookTorraNaru}). Thus, the functions performing this process (aggregation functions \cite{BookGraMariMesiPap,BookBeliPraCalvo}) have applications in decision theory, artificial intelligence \cite{KoProMi24}, economics \cite{GLMP15}, etc. In decision-making problems, aggregation functions are employed to aggregate individual preferences or opinions in the presence of uncertainty or conflicting information. For example, in multicriteria decision analysis (MCDA) (see, for example, \cite{LHSR22}), where decisions must be made based on several criteria, aggregation functions are used to merge different criteria weights and alternatives to determine the optimal solution. 

Easy examples of aggregation functions are the measures of central tendency in statistics such as the arithmetic mean, median, or mode. 

Additionally, we can also consider not only the aggregation of numbers but also the combination of more complex mathematical structures. In this context, Dobo\v{s} and his collaborators \cite{BorsikDobos81a,BorsikDobos81, Dobos95, Dobos98} explored the theory behind merging a family of metric spaces into a single metric space, where the ground set is the Cartesian product. To clarify this, let us consider a function $F:[0,+\infty)^I\to [0,+\infty).$ We say that $F$ is a \emph{metric preserving function} \cite{Dobos98} if for every family $\big\{(X_i,d_i):i\in I\big\}$ of metric spaces then $F\circ d_\Pi$ is a metric on $\prod_{i\in I} X_i,$ where $d_\Pi:\prod_{i\in I} X_i\times \prod_{i\in I} X_i\to [0,+\infty) ^I$ is given by
	$$d_\Pi(x,y)=(d_i(x_i,y_i))_{i\in I}$$
	for all $x=(x_i)_{i\in I}, y=(y_i)_{i\in I}\in\prod_{i\in I} X_i.$ The metric preserving functions were characterized and deeply studied in \cite{BorsikDobos81a,BorsikDobos81} (refer also to the surveys \cite{Corazza99,Dobos95}).

Similarly, Mayor and Valero \cite{MayorVale19} examined a related issue by characterizing the functions that aggregate multiple metrics defined on the same set into a single metric on that set. Concretely, they characterized those functions $F:[0,+\infty)^I\to [0,+\infty)$ verifying that whenever $\big\{(X,d_i):i\in I\big\}$ is a family of metric spaces then $F\circ d_\Delta$ is a metric on $X$ where $d_\Delta: X\times X\to [0,+\infty)^I$ is given by
	$$d_\Delta (x,y)=(d_i(x,y))_{i\in I}$$
	for all $x,y\in X.$

These two processes are distinct, and we differentiate them by referring to the first as aggregation on products and the second as aggregation on sets.

The aggregation problem is not exclusive to metrics. It has also been explored for other mathematical structures, such as quasi-metrics \cite{MayorVale10,MinyaVale19}, fuzzy quasi-pseudometrics \cite{PRLVale21}, norms \cite{HerMos91,PRL21}, asymmetric norms \cite{MartinMayorVale11,PRL21} or probabilistic quasi-uniformities \cite{PRL20b}, among others.

On the other hand, metric modulars were introduced by Chistyiakov \cite{Chis06,Chis10a,Chis10b,Chis15} as a gene\-ralization of Nakano's modulars to arbitrary sets. Roughly speaking, a metric modular is a metric that depends on a parameter $t\in (0,+\infty)$ (see Definition \ref{def:mqpm}). A typical example extracted from \cite{Chis15} is as follows: Given a metric space $(X,d)$ and $t>0,$ consider $w(t,x,y):=d(x,y)/t,$ which can be interpreted as the mean velocity between the points $x$ and $y$ over time $t.$ The function $w$ is the prototypical example of a metric modular and, by axiomatizing its fundamental properties, the definition of a metric modular emerged. Chystiakov also developed topological and convergence properties of metric modular spaces, demonstrating their coherence with the classical theory of modular linear spaces. Consequently, metric modular serves as 
	an important tool in nonlinear analysis. Furthermore, metric modulars have been studied without the symmetry axiom \cite{SebogodiPhD} and have been applied in fixed-point theory \cite{AlaEgePark16,DehEGEba12,EgeAla15,MongSintuKuman11}. Additionally,  a relationship with fuzzy metrics has also been established \cite{Roma24}. These studies highlight the significance of metric modulars across diverse areas.

We notice that the parameter included in metric modulars improves their flexibility,  making them more suitable for applications than classical metrics. For instance, consider a clustering problem where the objects to be classified are defined on different scales, that is, the data points are measured on different units \cite{HanKamberPei12}. In such cases, it is not appropriate to use a single measure to determine the proximity between the points. Adaptative and asymmetric distances, such as quasi-pseudometric modulars, can be particularly useful in this context (see \cite{Oza83} ). Motivated by this issue and applications in multi-agent systems, recent research by  Bibiloni-Femenias, Mi\~nana, and Valero has analyzed the problem of the aggregation of quasi-pseudometric modulars on sets in two papers \cite{BFMinyaVale23,BFVale24}. They characterized these functions and         showed that quasi-pseudometric modular aggregation functions on sets are coincident with the pseudometric modular aggregation functions on sets \cite[Theorem 8]{BFVale24}. This coincidence also occurs with the quasi-metric modular aggregation functions on sets and the metric modular aggregation functions on sets \cite[Theorem 9]{BFVale24}. However, the pro\-blem of characterizing quasi-pseudometric modular aggregation functions on products has not been addressed in the literature. The goal of this paper is to fill this gap by characterizing these functions (see Theorems \ref{thm:FpresifFnablapres}, \ref{thm:pmod_long} and \ref{thm:sepp_qmmodp}).
	Specifically, we demonstrate that the (quasi-)pseudometric modular aggregation functions on products mirror those on sets (Theorems \ref{thm:FpresifFnablapres} and \ref{thm:pmod_long}). Nevertheless, in the case of aggregating quasi-metric mo\-dulars, the two families differ (Example \ref{ex:sets_not_products}). In all cases, we clearly characterize these functions in terms of isotonicity and subadditivity. This leads to new insights in the theory of aggregation functions.
Our approach leverages the general framework of aggregation function theory developed in \cite{FBRL}. It is based on the fact that numerous mathematical structures whose aggregation functions have been characterized in the literature are indeed enriched categories over quantales (refer to Sections \ref{sec:quantales},\ref{sec:Qcategories}). In this way, in \cite{FBRL} the authors show that lax morphisms of quantales are an appropriate extension of the notion of aggregation functions and demonstrate that some results about the aggregation of metrics and fuzzy metrics can be inferred from this theory. Additionally, in \cite{LPPRL}, it is proven that quasi-pseudometric modular spaces are categorically isomorphic to enriched categories over the quantale $\nabla$ of nonincreasing functions $f:(0,+\infty)\to [0+\infty]$ (see Theorem \ref{thm:iso}). Consequently, taking advantage of the general theory of aggregation functions for quantales, we will characterize the quasi-pseudometric modular aggregation functions on products. Moreover, we will show that some of the results of \cite{BFMinyaVale23,BFVale24} follow from this theory.

The summary of the paper is as follows. In Section \ref{sec:quantales}, we compile the basic theory about quantales along with key examples. Section \ref{sec:Qcategories} addresses the core ideas of categories that are enriched over a quantale. It is highlighted that extended quasi-pseudometric spaces and fuzzy quasi-pseudometric spaces are forms of such enriched categories. In Section \ref{sec:lax}, we summarize some of the results of \cite{FBRL}, showing that lax morphisms of quantales are suitable functions for aggregating enriched categories over quantales. These results will be crucial to the aim of the paper. Section \ref{sec:qpmmod} introduces the mathematical structures we aim to aggregate, known as quasi-pseudometric mo\-dulars. These were originally defined by Chistyakov in \cite{Chis10a} to extend Nakano’s modular concept to arbitrary sets. Furthermore, we incorporate results from \cite{LPPRL}, which indicate that quasi-pseudometric modular spaces can be viewed as categories enriched over a specific quantale. With all this theory, we will characterize in Sections \ref{sec:aqpmod}, \ref{sec:apmod}, \ref{sec:aqmod} the (quasi-)(pseudo)metric modular aggregation functions on products.

%
\section{Quantales and lax morphisms}\label{sec:quantales}
%

In this section, we review the fundamental theory of quantales. Our primary references are \cite{BookQuantales,BookMonoidalTopology}.

\begin{definition}[\mbox{\cite[Section II.1.10]{BookMonoidalTopology},\cite[Section 2.3]{BookQuantales}}]
	A \textbf{quantale} $(\mathscr{V},\preceq,\ast)$ is a complete lattice $(\mathscr{V},\preceq)$ such that $\ast:\mathscr{V}\times \mathscr{V}\to \mathscr{V}$ is an associative binary operation which distributes over suprema:
	\begin{align*}
		u\ast \left(\bigvee_{i\in I} v_i\right)&=\bigvee_{i\in I}(u\ast v_i),\\
		\left(\bigvee_{i\in I} v_i\right)\ast u&=\bigvee_{i\in I}(v_i\ast u).
	\end{align*}
	If $\ast$ is also commutative then $(\mathscr{V},\preceq,\ast)$ is a \textbf{commutative quantale}.
	
	A quantale is called \textbf{unital} if $\ast$ has a unit $1_\mathscr{V}.$ A unital quantale is \textbf{integral} if the unit is the top element $\top$ of $\mathscr{V}.$
	
\end{definition}

In the remainder of the paper, we will only consider commutative integral quantales; however, for simplicity, we will use only the term quantale.

\begin{remark}\label{rem:1}
	Notice that in an integral quantale $(\mathscr{V},\preceq,\ast)$ we have that $u\ast v\preceq u\wedge v$ for all $u,v\in \mathscr{V}.$ In fact, $u=u\ast (v\vee 1_\mathscr{V})=(u\ast v)\vee (u\ast 1_\mathscr{V})=(u\ast v)\vee u$ so $u\ast v\preceq u.$ In a similar way, $u\ast v\preceq v.$
\end{remark}

\begin{example}[\mbox{\cite[Example II.1.10.1]{BookMonoidalTopology}}]\label{ex:quantales}\
	\begin{enumerate}
		\item Let $\ast$ be a triangular norm (t-norm for short) on $[0,1]$ (see, for example, \cite{BookTN}), that is, an associative, commutative binary operation $\ast\colon[0,1]\times [0,1]\rightarrow [0,1]$ with unit $1$, such that $a\ast b\leq c\ast d$ whenever $a\leq c$ and $b\leq d$, with $a,b,c,d\in [0,1]$. If $\ast$ is left-continuous then $\big([0,1],\leq,\ast\big)$ is a commutative integral quantale, where $\leq$ is the usual order.
		\item Let $\mathbf{2}=\{0,1\}$ be a set with two different elements endowed with the usual order. If $\ast$ is an arbitrary t-norm then $(\mathbf{2},\leq,\ast)$ is a commutative integral quantale.
		\item  Let us consider the opposite order $\leq^\text{\normalfont op}$ on the extended real line $[0,+\infty].$ Specifically, $x\leq^\text{\normalfont op} y$ if and only if $y\leq x$. If we extend the usual sum $+$ on the real numbers to include $+\infty$ as usual, then $\mathsf{P}_+=\big([0,+\infty],\leq^\text{\normalfont op},+\big)$ forms a commutative integral quantale known as the \textbf{Lawvere quantale} \cite{CookWeiss22} (see also \cite[Example II.1.10.1.(3)]{BookMonoidalTopology}).


		\item Let $\Delta_+$ be the family of distance distribution functions given by
		$$\Delta_+=\Big\{f:[0,+\infty]\to [0,1]: f\text{ is isotone and left-continuous} \Big\},$$
		where left-continuous means that $f(x)=\bigvee_{y<x}f(y)$ for all $x\in [0,+\infty]$ (as usual, $\bigvee\varnothing =0$). Then $\Delta_+,$ endowed with the pointwise order $\leq,$ is a complete lattice. 
		Moreover, given a left-continuous t-norm $\ast$, $f,g\in\Delta_+ $ and $t\in [0,\infty]$ define $f\circledast g\in \Delta_+ $ as
		$$(f\circledast g)(t)=\bigvee_{r+s\leq t}f(r)\ast g(s)=\bigvee_{r+s= t}f(r)\ast g(t).$$
		Then $(\Delta_+,\leq,\circledast)$ is a quantale (see \cite{Flagg97,HofReis13,JagerShi19}) where the unit $f_{0,1}:[0,+\infty]\to [0,1]$ is given by
		$$f_{0,1}(t)=\begin{cases}
			0&\text{ if }t=0,\\
			1&\text{ otherwise},
		\end{cases}$$
		for all $t\in [0,+\infty].$ For simplicity, we will denote the quantale $(\Delta_+,\leq,\circledast)$ by $\Delta_+(\ast).$
		\item 	Let us consider the set 
		$$\nabla:=\Big\lbrace f:(0,+\infty) \to [0,+\infty]^{op} \text{ : } f \text{ isotone }\Big\rbrace$$
		endowed with the pointwise order also denoted by $\leq^{op}$. Then 	$(\nabla,\leq^{op},\oplus)$ is a quantale \cite{LPPRL}, where $\oplus : \nabla \times \nabla \to \nabla$ is given by:
		$$(f\oplus g)(t):=\bigvee_{r+s\leq t}^{\leq^{op}}(f(r)+g(s))=\bigvee_{r+s= t}^{\leq^{op}}(f(r)+g(s)).$$
		Its unit is the constant 0 function that we denote by $0_\nabla.$ 
		
		This quantale will be crucial in the paper because, as proved in \cite{LPPRL}, quasi-pseudometric modular spaces can be viewed as categories enriched over $\nabla$ (see Theorem \ref{thm:iso}).
	\end{enumerate}
\end{example}

\begin{remark}\label{rem:product_quantales}
	Let $\Big\{(\mathscr{V}_i,\preceq_i,\ast_i):i\in I\Big\}$ be an arbitrary family of quantales. Define $\preceq:=\Pi_{i\in I}\preceq_i$ as the product partial order defined componentwisely, and let $\ast:=\Pi_{i\in I}~\ast_i :(\prod_{i\in I}\mathscr{V}_i)\times(\prod_{i\in I}\mathscr{V}_i)\to \prod_{i\in I}\mathscr{V}_i$ be the componentwise operation given by $(u\ast v)_i=u_i\ast_i v_i$ for all $u=(u_i)_{i\in I}, v=(v_i)_{i\in I}\in \prod_{i\in I}\mathscr{V}_i$ and all $i\in I.$ Thus, $(\prod_{i\in I}\mathscr{V}_i,\preceq,\ast)$ forms a quantale. 
	
	If all the quantales in the family $\Big\{(\mathscr{V}_i,\preceq_i,\ast_i):i\in I\Big\}$ are equal to $(\mathscr{V},\preceq,\ast),$ we will simplify the notation by also using $\preceq$ and $\ast$ to denote the partial order and the operation on $\mathscr{V}^I.$ 
	
	In particular, we denote by $\mathsf{P}_+^I$ the quantale $\big([0,+\infty]^I,\leq^\text{\normalfont op},+\big)$ and by $\nabla^I$ the quantale $(\nabla^I,\leq^{\text{op}},\oplus).$ The unit of $\mathsf{P}_+^I$ will be denoted by $0_I,$ which re\-presents $(0)_{i\in I}.$ Furthermore, $0_{\nabla^I}$ designates the unit of $\nabla^I$, that is, $(0_\nabla)_{i\in I}.$
\end{remark}

By considering quantales as ordered categories, the concept of lax functor reduces to the following:

\begin{definition}[\mbox{\cite{HofReis13},\cite[Section II.4.6]{BookMonoidalTopology}}]
	A map $F:(\mathscr{V},\preceq,\ast)\to(\mathscr{W},\curlyeqprec,\star)$ between two quantales is said to be a \emph{lax morphism of quantales} if
	\begin{itemize}
		\item $u\preceq v$ implies $F(u)\curlyeqprec F(v)$ for all $u,v\in\mathscr{V};$\hfill (isotone)
		\item $1_\mathscr{W}\curlyeqprec F(1_\mathscr{V});$
		\item $F(u)\star F(v)\curlyeqprec F(u\ast v)$ for all $u,v\in\mathscr{V}.$\hfill (subadditive)
	\end{itemize}
	
\end{definition}

\begin{example}\label{ex:lax_morphism_P}
	If $\mathscr{V}=\mathsf{P}_+^I$ and $\mathscr{W}=\mathsf{P}_+$, then a function $F:\mathsf{P}_+^I\to\mathsf{P}_+$ is a lax morphism if 
	\begin{itemize}
		\item $F$ is isotone;
		\item $F(0_I)=0$;
		\item $F(x+y)\leq F(x)+F(y)$, for all $x,y\in [0,+\infty]^I$ ($F$ is subadditive).
	\end{itemize}
	Notice that, in the last inequality, the symbol $+$ denotes the usual sum of extended real numbers on the right side, while on the left side, it indicates the sum performed componentwise.
\end{example}

\section{Categories enriched over a quantale}\label{sec:Qcategories}

As previously discussed, numerous mathematical structures for which the aggregation functions have been characterized serve as examples of categories enriched over a quantale. This section aims to gather the foundational concepts of this theory, beginning with the definition of a $\mathscr{V}$-category, which is a category enriched over a quantale $\mathscr{V}$.

\begin{definition}[\mbox{\cite[Section III.1.3]{BookMonoidalTopology}, c.f. \cite[Definition 3.1]{FlaggKopp97}}]
	Let $(\mathscr{V},\preceq,\ast)$ be a quantale. A \emph{$\mathscr{V}$-category} is a pair $(X,a)$ where $X$ is a nonempty set and $a:X\times X\to\mathscr{V}$ is a map such that
	\begin{itemize}
		\item[(VC1)]\; $1_\mathscr{V}\preceq a(x,x)$ for all $x\in X;$
		\item[(VC2)]\; $a(x,z)\ast a(z,y)\preceq a(x,y)$ for all $x,y,z\in X.$
	\end{itemize}
	Moreover, a $\mathscr{V}$-category $(X,a)$ is said to be:
	\begin{itemize}
		\item \emph{separated} if given $x,y \in X$, whenever $1_\mathscr{V}\preceq a(x,y)$ and $1_\mathscr{V}\preceq a(y,x)$ then $x=y.$
		\item \emph{symmetric} if $a(x,y)=a(y,x)$ for all $x,y\in X.$
	\end{itemize}
	
	A function $f:(X,a)\to (Y,b)$ between two $\mathscr{V}$-categories is called a \emph{$\mathscr{V}$-functor} if
	$$a(x,y)\preceq b(f(x),f(y))$$
	for all $x,y\in X.$
	
	The category whose objects are $\mathscr{V}$-categories and whose morphisms are $\mathscr{V}$-functors will be denoted by $\mathscr{V}$-$\mathsf{Cat}.$
	
\end{definition}

\begin{remark}
	It is noteworthy that each topological space can be interpreted as a $\mathscr{V}$-category (refer to \cite{Flagg97,Kopp88}).
\end{remark}

We will now provide several key examples of $\mathscr{V}$-categories to illustrate that they encompass some important mathematical structures.

\begin{example}[\mbox{\cite[Examples III.1.3.1]{BookMonoidalTopology},\cite{HofReis13}}]\label{ex:vcat}\
	\begin{enumerate}
		\item $\mathbf{2}$-categories can be interpreted as preordered sets, while $\mathbf{2}$-functors co\-rrespond to isotone (preserving order) maps. Specifically, for a $\mathbf{2}$-category $(X,a)$, define the binary relation $\preceq_a$ on $X$ given by $x\preceq_a y$ if and only if $a(x,y)=1,$ for all $x,y\in X.$ This relation is a preorder on $X.$ 
		
		Conversely, if $(X,\preceq)$ is a partially ordered set then $a_\preceq:X\times X\to\mathscr{V}$ given by $a_{\preceq}(x,y)=1$ if $x\preceq y$ and $a_{\preceq}(x,y)=0$ otherwise, turns $(X,\preceq)$ into a $\mathbf{2}$-category $(X,a_\preceq)$.
		As a result, $\mathbf{2}$-$\mathsf{Cat}$ is isomorphic to the category $\mathsf{POrd}$ of preordered sets. 
		
		Obviously, the category of separated $\mathbf{2}$-categories is isomorphic to the category of partially ordered sets.

		\item  $\mathsf{P}_+$-categories are equivalent to extended quasi-pseudometric spaces and $\mathsf{P}_+$-morphisms are equivalent to non-expansive maps, where $\mathsf{P}_+$ is the Lawvere quantale (see Example \ref{ex:quantales}).
		
		\noindent In fact, the axioms of a $\mathsf{P}_+$-category can be expressed as follows:
		\begin{quote}
			\begin{itemize}
				\item[(VC1)]\; $ 0\leq^{op} a(x,x)\Leftrightarrow  0\geq a(x,x)\Leftrightarrow 0=a(x,x);$
				\item[(VC2)]\; $a(x,y)+a(y,z)\leq^{op}a(x,z)\Leftrightarrow a(x,z)\leq a(x,y)+a(y,z),$
			\end{itemize}
		\end{quote}
		for all $x,y,z\in X$. Thus, $a$ is an extended quasi-pseudometric on $X$ \cite{Beer13,OtaTokoMuko19}. The term ``extended'' indicates that $a$ may take the value $+\infty$; ``quasi'' emphasizes the fact that $a$ is not necessarily symmetric; and ``pseudo'' means that $a$ does not verify the separated property of a $\mathscr{V}$-category.
		
		It is also clear that an extended quasi-pseudometric is a $\mathsf{P}_+$-category. Therefore, $\mathsf{P}_+$-categories and extended quasi-pseudometric spaces are equivalent concepts. 
		
		In a similar vein, separated $\mathsf{P}_+$-categories and symmetric separated $\mathsf{P}_+$-categories are equivalent to extended quasi-metric spaces and extended metric spaces respectively.
		
		Moreover, a $\mathsf{P}_+$-functor $f:(X,a)\to (Y,b)$ verifies
		$$a(x,y)\geq b(f(x),f(y))$$
		for all $x,y\in X$, that is, $f$ is a non-expansive mapping. Hence, $\mathsf{P}_+$-$\mathsf{Cat}$ is isomorphic to the category $\mathsf{QPMet}$ of quasi-pseudometric spaces. 
		
		\item If we consider the quantale $\Delta_+(\ast):=(\Delta_+,\leq,\circledast
		)$ (see Example \ref{ex:quantales}.(4)), where $\ast$ is a continuous t-norm, then the category $\Delta_+(\ast)$-$\mathsf{Cat}$ is isomorphic to the category $\mathsf{FQPMet}(\ast)$ of fuzzy quasi-pseudometric spaces \cite{GregRo04} with respect to $\ast$ and fuzzy nonexpansive maps. 
		
		Notice that if $(X,a)$ is a $\Delta_+(\ast)$-category then $M_a:X\times X\times [0,+\infty]
		\to [0,1]$ defined as $M_a(x,y,t)=a(x,y)(t)$ for every $x,y\in X,$ $t\in [0,+\infty]$, that is, the evaluation of $a(x,y)$ at $t$, satisfies:
		\begin{quote}
			\begin{itemize}
				\item[(FM1)]\; $M_a(x,y,0)=a(x,y)(0)=0$ for all $x,y\in X;$
				\item[(FM2)]\; $M_a(x,x,t)=1$ for all $t>0;$
				\item[(FM3)]\; $M_a(x,z,s)\ast M_a(z,y,t)\leq M(x,y,s+t)$ for every $x,y,z\in X$, $s,t\in [0,+\infty];$
				\item[(FM4)]\; $M_a(x,y,\cdot)=a(x,y)$ is left-continuous for every $x,y\in X.$
			\end{itemize}
		\end{quote}
		The previous axioms are exactly the conditions that a pair $(M,\ast)$ must satisfy for being a fuzzy quasi-pseudometric on $X$ (see, for example, \cite{GregRo04}), where $M:X\times X\times [0,+\infty)\to [0,1]$ and $\ast$ is a continuos t-norm. Observe that $M_a(x,y,t)$ also exists when $t=\infty$ and $M_a(x,y,\infty)=\bigvee_{0\leq t<\infty} M_a(x,y,t).$ This is not required in the classical definition of a fuzzy quasi-pseudometric but there is no loss of generality if we add this property.\\
		Conversely, let $(M,\ast)$ be a fuzzy quasi-pseudometric on $X.$ If we define $m:X\times X\to\Delta_+(\ast)$ as
		$$m(x,y)(t)=\begin{cases}
			M(x,y,t) &\text{ if } 0\leq t < +\infty,\\
			\bigvee\limits_{0\leq s<+\infty} M(x,y,s)&\text{ if } t=+\infty,
		\end{cases}$$
		for every $x,y\in X$, $t\in [0,+\infty]$, then it is straightfoward to check that $(X,m)$ is a $\Delta_+(\ast)$-category. Observe that $m$ is the exponential mate $\ceil{M}$ of the extension of $M$ to $X\times X\times [0,+\infty].$
		
	\end{enumerate}
\end{example}

		
		


\begin{example}\label{exmp:prodq}\
	Let $(\mathscr{V},\preceq,\ast)$ be a quantale.
	\begin{enumerate}
		\item  If $(X,a)$ is a $\mathscr{V}^I$-category, for each $i\in I$ let $a_i:X\times X\to \mathscr{V}$ be the $i$th-coordinate function of $a$, that is, $a_i(x,y)=(a(x,y))_i$ for every $x,y\in X.$ Then $\big\{(X,a_i):i\in I\big\}$ is a family of $\mathscr{V}$-categories.
		If $(X,a)$ is a symmetric $\mathscr{V}^I$-category then $(X,a_i)$ is a symmetric $\mathscr{V}$-category for all $i\in I.$
		\item  Let $\big\{(X_i,a_i):i\in I\big\}$ be a family of $\mathscr{V}$-categories. Define $a_\Pi:\prod_{i\in I} X_i\times \prod_{i\in I} X_i\to \mathscr{V}^I$ as
		$$(a_\Pi(x,y))_i=a_i(x_i,y_i)$$
		for all $x,y\in \prod_{i\in I} X_i, i\in I.$ It is straightforward to check that $(\prod_{i\in I} X_i,a_\Pi)$ is a $\mathscr{V}^I$-category.
		\item Let $\big\{(X,a_i):i\in I\big\}$ be a family of $\mathscr{V}$-categories. Define $a_\Delta: X\times X\to \mathscr{V}^I$ as
		$$(a_\Delta(x,y))_i=a_i(x,y)$$
		for all $x,y\in X, i\in I.$ Then $(X,a_\Delta)$ is a $\mathscr{V}^I$-category.
	\end{enumerate}
\end{example}

%
\section{Lax morphisms of quantales as generalized aggregation functions}\label{sec:lax}
%

In \cite{FBRL}, the authors expand upon the existing theory of aggregation functions by situating it within the framework of categories enriched over a quantale. This conceptual shift allows for a more robust analysis of aggregation functions for various mathematical structures.

In our paper, we aim to leverage this generalized framework to explore and derive results concerning the aggregation of quasi-pseudometric modulars (see Section \ref{sec:qpmmod} for a more in-depth discussion of this concept).

To ensure a thorough understanding of our subsequent findings, we will first review the results from \cite{FBRL} that will be instrumental in our research later in the paper. 

	
	
	

\begin{definition}[\cite{FBRL}]\label{def:preserving}
	A map $F:(\mathscr{V},\preceq,\ast)\to(\mathscr{W},\curlyeqprec,\star)$ between two quantales is said to be \textbf{preserving} if the map $\mathsf{F}:\mathrm{Obj}(\mathscr{V}\text{-}\mathsf{Cat})\to\mathrm{Obj}(\mathscr{W}\text{-}\mathsf{Cat}),$ which assigns to a $\mathscr{V}$-category $(X,a)$ the pair $(X,F\circ a)$ is well-defined, that is, if $(X,F\circ a)$ is a $\mathscr{W}$-category.
	
	We will denote by $\mathscr{P}\big((\mathscr{V},\preceq,\ast),(\mathscr{W},\curlyeqprec,\star)\big)$, or simply by $\mathscr{P}(\mathscr{V},\mathscr{W})$ if no confusion arises, the family of preserving functions between the quantales $(\mathscr{V},\preceq,\ast)$ and $(\mathscr{W},\curlyeqprec,\star).$
	
	If the map $F$ satisfies that $(X,F\circ a)$ is a separated (resp. symmetric) $\mathscr{W}$-category whenever $(X,a)$ is a separated (resp. symmetric) $\mathscr{V}$-category then $F$ is said to be \textbf{separately preserving} (resp. \textbf{symmetrically preserving}).  The family of separately (resp. symmetrically) preserving functions will be denoted by $\mathsf{se}\mathscr{P}\big((\mathscr{V},\preceq,\ast),(\mathscr{W},\curlyeqprec,\star)\big)$ (resp. $\mathsf{sy}\mathscr{P}\big((\mathscr{V},\preceq,\ast),(\mathscr{W},\curlyeqprec,\star)\big)$).
	
	Notice that $\mathscr{P}(\mathscr{V},\mathscr{W})\subseteq \mathsf{sy}\mathscr{P}(\mathscr{V},\mathscr{W})$.
\end{definition}

We now present an easy example of the family of preserving functions $\mathscr{P}(\mathscr{V},\mathscr{W})$ for specific quantales $\mathscr{V}$ and $\mathscr{W}$. Further examples will be provided in Theorems \ref{thm:FpresifFnablapres}, \ref{thm:pmod_long}, \ref{thm:sepp_qmmodp}, and Proposition \ref{prop:preservingnabla} 
(see also \cite{FBRL}).

\begin{example}
	Let $F:(\mathbf{2},\leq,\ast)\to (\mathbf{2},\leq,\ast)$ be a function. Then $F$ is preserving if and only if $F$ is the identity or $F$ is the identically 1 function. Otherwise, $F(1)=0$ and in this case, given a $\mathbf{2}$-category $(X,a)$ (a partially ordered set, Example \ref{ex:vcat})  and $x\in X$, then $(X,F\circ a)$ is not a $\mathbf{2}$-category since $1\not\leq F(a(x,x))=F(1)=0$ failing to fulfill (VC1).
\end{example}

The concept of triangle triplet, which first appeared in \cite{Terpe84}, plays a crucial role in the characterization of quasi-pseudometric aggregation functions \cite{Dobos98,MayorVale10}. A counterpart notion appears in \cite{PRLVale21} for characterizing fuzzy quasi-pseudometric aggregation functions. The following concept expands this notion to a more general context.

\begin{definition}[\cite{FBRL}]
	Let $(S,\leq,\cdot)$ be an ordered semigroup, that is, a semigroup endowed with a partial order compatible with the operation. A triplet $(x,y,z)\in S^3$ is said to be an \textbf{asymmetric triangle triplet} on $(S,\leq,\cdot)$ if $x\cdot y\leq z.$ Moreover, it is said to be a \textbf{triangle triplet} if every permutation of the triplet is an asymmetric triangle triplet.
	
	A function $F:(S,\leq,\cdot)\to(M,\leqslant,\ast)$ between two ordered semigroups is said to \textbf{preserve (asymmetric) triangle triplets} if $(F(x),F(y),F(z))$ is a(n) (asymmetric) triangle triplet on $(M,\leqslant,\ast)$ whenever $(x,y,z)$  is a(n) (asymmetric) triangle triplet on $(S,\leq,\cdot).$
\end{definition}

\begin{example}[\mbox{\cite[Chapter 2.2]{Dobos98}}]\
	\begin{itemize}
		\item A triangle triplet on $\mathsf{P}_+$ is a triplet $(a,b,c)\in [0,+\infty]^3$ such that 
		$$a\leq b+c,\hspace*{0.5cm}b\leq a+c,\hspace*{0.5cm}c\leq a+b.$$
		An asymmetric triangle triplet on $\mathsf{P}_+$ is a triplet $(a,b,c)\in [0,+\infty]^3$ such that 
		$$a\leq b+c.$$
		
		\item if $(X,d)$ is a (quasi-)pseudometric space, then $\big(d(x,y),d(y,z),$ $d(x,z)\big)$ is a(n) (asymmetric) triangle triplet on $\mathsf{P}_+$ for every $x,y,z\in X.$   
	\end{itemize}
\end{example}

\begin{example}
	Every constant function $F:(\mathscr{V},\preceq,\ast)\to (\mathscr{W},\curlyeqprec,\star)$ between two commutative integral quantales preserves asymmetric triangle triplets. Notice that whenever $(x,y,z)\in\mathscr{V}^3$ is an asymmetric triangle triplet, then $(F(x),F(y),F(z))=(F(x),F(x),F(x))$ is a triangle triplet on $\mathscr{W}$ (see Remark \ref{rem:1}).
\end{example}

The results detailed below, taken from \cite{FBRL}, characterize the various families of preserving functions that were introduced earlier.

\begin{theorem}[\mbox{\cite{FBRL}}]\label{thm:plmq}
	Let $(\mathscr{V},\preceq,\ast),(\mathscr{W},\curlyeqprec,\star)$ be two commutative integral quantales. The following statements are equivalent:
	\begin{enumerate}
		
		\item  $F:(\mathscr{V},\preceq,\ast)\to(\mathscr{W},\curlyeqprec,\star)$ is preserving; 
		\item $F:(\mathscr{V},\preceq,\ast)\to(\mathscr{W},\curlyeqprec,\star)$ is a lax morphism;
		\item $F$ preserves asymmetric triangle triplets and $1_\mathscr{W}=F(1_\mathscr{V})$;
		\item $F$ preserves symmetric triangle triplets, it is isotone and $1_\mathscr{W}=F(1_\mathscr{V}).$
	\end{enumerate}
	
\end{theorem}

\begin{theorem}[\mbox{\cite{FBRL}}]\label{thm:symp}
	Let $(\mathscr{V},\preceq,\ast),(\mathscr{W},\curlyeqprec,\star)$ be two commutative integral quantales. Then a function $F:(\mathscr{V},\preceq,\ast)\to(\mathscr{W},\curlyeqprec,\star)$ is symmetrically preserving if and only if $F$ preserves triangle triplets and $1_\mathscr{W}=F(1_\mathscr{V}).$
\end{theorem}

\begin{proposition}[\mbox{\cite{FBRL}}]\label{prop:asym_symiso}
	Let $(\mathscr{V},\preceq,\ast),(\mathscr{W},\curlyeqprec,\star)$ be two commutative integral quantales. Then $F$ preserves asymmetric triangle triplets and $1_\mathscr{W}=F(1_\mathscr{V})$ if and only if $F$ preserves symmetric triangle triplets, is isotone and $1_\mathscr{W}=F(1_\mathscr{V}).$
\end{proposition}

\begin{theorem}[\mbox{\cite{FBRL}}]\label{thm:seplmq}
	
	Let $(\mathscr{V},\preceq,\ast),(\mathscr{W},\curlyeqprec,\star)$ be two commutative integral quantales. The following statements are equivalent:
	\begin{enumerate}
		\item  $F:(\mathscr{V},\preceq,\ast)\to(\mathscr{W},\curlyeqprec,\star)$ is separately preserving; 
		\item $F:(\mathscr{V},\preceq,\ast)\to(\mathscr{W},\curlyeqprec,\star)$ is a lax morphism satisfying $F^{-1}\left(1_{\mathscr{W}}\right)=\{1_{\mathscr{V}}\}$;
		\item $F$ preserves asymmetric triangle triplets and $F^{-1}\left(1_{\mathscr{W}}\right)=\{1_{\mathscr{V}}\}.$
	\end{enumerate}
\end{theorem}

%
\section{Quasi-pseudometric modular spaces}\label{sec:qpmmod}
%

Recall that our objective is to characterize the aggregation functions of a mathematical structure known as quasi-pseudometric modular. In this section, we review this concept originally defined by Chistyiakov \cite{Chis06,Chis10a,Chis10b,Chis15}, as well as its asymmetric version, which was explored in \cite{SebogodiPhD}. To obtain the promised characterizations, we will draw upon the general theory of aggregation functions developed in the previous results. However, this requires us to describe quasi-pseudometric modular spaces as categories enriched over a specific quantale. This topic has been addressed in \cite{LPPRL}, from which we will gather relevant results.

\begin{definition}[\cite{Chis15,SebogodiPhD}]\label{def:mqpm}
	Let $X$ be a non-empty set. A function
	$$w:(0,+\infty)\times X\times X \to [0,+\infty]$$
	is said to be a \textbf{quasi-pseudometric modular} on $X$ if for every $x,y,z \in X$ and $t,s >0$ it verifies:
	\begin{enumerate}
		\item[(M1)] $w(t,x,x)=0$ for every $t>0$,
		\item[(M2)] $w(t+s,x,y)\leq w(t,x,z)+w(s,z,y)$.
	\end{enumerate}
	If, in addition, $w$ satisfies
	$$\text{(M3) }w(t,x,y)=w(t,y,x)=0 \text{ for every } t>0 \text{ if and only if } x=y,$$
	then $w$ is called a \textbf{quasi-metric modular}.\\
	If a quasi-(pseudo)metric modular $w$ verifies:
	$$\text{(M4) }w(t,x,y)=w(t,y,x) \text{ for every } x,y \in X \text{ and every }t>0,$$
	then $w$ is said to be a \textbf{(pseudo)metric modular} on $X$.\\
	
	The pair $(X,w)$ will be called a (quasi-)(pseudo)metric modular space.

\end{definition}

\begin{example}[\mbox{compare with \cite[Section 1.3.1]{Chis15}}]\label{exmp:wgd}
	Let $(X, d)$ be a quasi-pseudometric space, and $g:(0, +\infty) \rightarrow[0, +\infty]$ be a nonincreasing function. Define $w_{g,d}:(0,+\infty)\times X\times X\to [0,+\infty]$ as 
	$$w_{g,d}(t,x, y)=g(t) \cdot d(x, y),$$
	for all $t>0$, $x,y\in X,$ with the convention that $\infty \cdot 0=0$, and $\infty \cdot a=\infty$ for all $a>0$. Then $(X,w_{g,d})$ is a quasi-pseudometric modular space. 
	
	Moreover, if $g$ is not identically zero and $d$ is a quasi-metric, then $(X,w_{g,d})$ is a quasi-metric modular space. 
	
	When $g$ is the function identically 1, then
	$$w_{g,d}(t,x,y)=d(x,y)$$
	for all $t>0 0$, $x,y\in X.$ We denote this quasi-pseudometric modular by $w_d.$

\end{example}

\begin{example}\label{exmp:mods}
	Let $w:(0,+\infty)\times [0,+\infty)\times [0,+\infty)\to [0,+\infty]$ defined as
	$$w(t,x,y)=\begin{cases}
		x+y&\text{ if }x\neq y,\\
		0&\text{ if }x=y,
	\end{cases}$$
	for all $t>0,x,y\geq 0.$ Then it is easy to check that $([0,+\infty),w)$ is a metric modular space.
\end{example}

In \cite{LPPRL}, it has been shown that quasi-pseudometric modular spaces can be viewed as enriched categories over a quantale. We collect some of the results from \cite{LPPRL} in the following one.

\begin{theorem}[\cite{LPPRL}]\label{thm:iso}
	Let $\nabla$-$\mathsf{Cat}$ denote the family of all $\nabla$-categories meanwhile $\mathsf{QPMod}$ denotes the family of all quasi-pseudometric modular spaces. Define $\mathscr{I}:\nabla\text{-}\mathsf{Cat}\to \mathsf{QPMod}$ as 
	$$\mathscr{I}\big((X,a)\big)=(X,w_a)$$ for every $\nabla$-category $(X,a)$, where $w_a(t,x,y)=a(x,y)(t)$ for all $x,y\in X$, $t>0.$ 
	
	Then $\mathscr{I}$ is an isomorphism having as inverse
	$$\mathscr{I}^{-1}((X,w))=(X,a_w)$$
	where $a_w:X\times X\to \nabla$ is defined as
	$$a_w(x,y)(t)=w(t,x,y)$$
	for all $x,y\in X$, $t>0.$
\end{theorem}

\begin{remark}
	We note that the results proved in \cite{LPPRL} are more general because they consider the categories of quasi-pseudometric modular spaces and $\nabla$-categories rather than focusing solely on their objects. Nevertheless, in this study, we do not require that level of generality since we will not consider morphisms in our future discussion.
\end{remark}

\begin{remark}\label{rem:symsep}		
	A straightforward computation shows that indeed, separated $\nabla$-categories can be seen as quasi-metric modular spaces. Likewise, symmetric $\nabla$-categories can be seen as pseudometric modular spaces.
\end{remark}

\section{Aggregation functions of quasi-pseudometric modulars}\label{sec:aqpmod}
In this section and the following ones, we finally address the goal of this paper. As previously commented, in  \cite{BFMinyaVale23}, the concept of a (pseudo)metric modular aggregation function was introduced and characterized. This topic was further developed in \cite{BFVale24}, which explored quasi-pseudometric modular aggregation functions. The investigations conducted in \cite{BFMinyaVale23} and \cite{BFVale24} focus solely on what we refer to as aggregation functions on sets. This means that these functions enable the merging of a family of quasi-pseudometric modulars defined on the same set, into a single one in that set. 

Note that there is an alternative method for aggregating quasi-pseudometric modulars by constructing a new one within the Cartesian product. In this section, we will characterize these functions using the general theory of aggregation functions between quantales.
Moreover, we will prove some already known results of aggregation of (quasi-)(pseudo)metric modulars using the theory of quantales given in Section \ref{sec:lax}.
We will begin by clearly defining the two different concepts of functions that aggregate quasi-pseudometric modulars.

\begin{definition}[\mbox{compare with \cite{BFMinyaVale23,BFVale24}}]\label{def::ModOPModOS}
	A function $F:[0,+\infty]^{I} \to [0,+\infty]$ is said to be:
	\begin{itemize}
		\item a \textbf{(quasi-)(pseudo)metric modular aggregation function on products} if for every family $\big\lbrace (X_{i}, w_{i}) \text{ : } i \in I\big\rbrace$ of (quasi-)(pseudo)metric modular spaces, then $F \circ w_\Pi$ is a (quasi-)(pseudo)metric modular on $\Pi_{i \in I}X_{i}$, where
		$$w_\Pi:(0,+\infty)\times \prod_{i \in I}X_{i} \times \prod_{i \in I}X_{i}  \to [0,+\infty]^{I} $$ 
		is given by
		$$w_\Pi(t,x,y)=\big(w_{i}(t,x_{i},y_{i})\big)_{i \in I}$$
		for all $x,y\in\prod _{i\in I} X_i, t>0.$
		\item a \textbf{(quasi-)(pseudo)metric modular aggregation function on sets} if for every collection $\lbrace w_{i}: i \in I \rbrace$ of (quasi-)(pseudo)metric modulars over a nonempty set $X$, then $F \circ w_\Delta$ is a (quasi-)(pseudo)metric modular on $X$, where
		$$w_\Delta: (0,+\infty) \times X \times X  \to [0,+\infty]^{I}$$
		is given by
		$$w_\Delta(t,x,y)=\big(w_{i}(t,x,y)\big)_{i \in I}$$
		for all $x,y\in X, t>0.$
	\end{itemize}
We shall now give some examples that illustrate this definition.
	
	\begin{example}\
		
		\begin{enumerate}
			\item Given $k\in [0,+\infty],$ let $F_k:[0,+\infty]^I\to [0,+\infty]$ given as
				$$F_k(x)=\begin{cases}k&\text{ if }x\neq 0_i,\\
					0&\text{ if }x=0,\end{cases}$$
				for every $x\in [0,+\infty]^I.$ Then it is straightforward to check that $F_k$ is a quasi-pseudometric modular aggregation function on products.
			\item Given $n\in\mathbb{N},$ define $S_n:[0,+\infty]^n\to [0,+\infty]$ as
				$$S_n(x)=x_1+\ldots+x_n$$
				for all $x=(x_1,\ldots,x_n)\in [0,+\infty]^n.$ An easy computation shows that $S_n$ is a quasi-pseudometric modular aggregation function on products.
		\end{enumerate}
	\end{example}

	The family of quasi-pseudometric modular aggregation functions on products (resp. on sets) will be denoted by $\mathsf{QPModAP}$ (resp. $\mathsf{QPModAS}$). The notations $\mathsf{QModAP},$ $\mathsf{QModAS},$ $\mathsf{PModAP},$ $\mathsf{PModAS},$ $\mathsf{ModAP},$ $\mathsf{ModAS}$ are self-explained.
	
\end{definition}
The families $\mathsf{PModAS}, \mathsf{ModAS}$ were characterized in \cite{BFMinyaVale23} meanwhile the families $\mathsf{QPModAS}, \mathsf{QModAS}$ were characterized in \cite{BFVale24}. In this paper, we will characterize the families $\mathsf{QPModAP}, \mathsf{QModAP}, \mathsf{PModAP}, \mathsf{ModAP}$. To achieve this, we will make use of the general theory about aggregation functions developed by the authors in \cite{FBRL}. Moreover, we will demonstrate that some of the results of \cite{BFMinyaVale23,BFVale24} can be deduced from that general theory.  

We start by presenting a result that demonstrates the equivalence between the quasi-pseudometric modular aggregation functions on products and on sets. Additionally, it establishes a helpful connection between the classical theory of aggregation of quasi-pseudometric modulars and the theory of aggregation for $\nabla$-categories, showing that a quasi-pseudometric modular aggregation function on products $F$ induces a function transforming $\nabla^I$-categories into $\nabla$-categories.

\begin{proposition}\label{prop:cons_f_nabla}
	Let $F:[0,+\infty]^I\to [0,+\infty]$ be a function. Define $F_\nabla: \left([0,+\infty]^{(0,+\infty)}\right)^I\to [0,+\infty]^{(0,+\infty)}$ as
	$$F_\nabla\big((f_i)_{i\in I}\big)(t)=F\big((f_i(t))_{i\in I}\big)$$
	for every $(f_i)_{i\in I}\in \left([0,+\infty]^{(0,+\infty)}\right)^I$ and $t>0.$ Then the following statements are equivalent:
	\begin{enumerate}
		\item $F\in \mathsf{QPModAP};$
		\item $F\in \mathsf{QPModAS};$
		\item whenever $(X,a)$ is a $\nabla^I$-category then $(X,F_\nabla\circ a)$ is a $\nabla$-category.
	\end{enumerate}
	
\end{proposition}

\begin{proof}
	(1) $\Rightarrow$ (2)
	This is straightforward. 
	
	(2) $\Rightarrow$ (3) Suppose that $F\in \mathsf{QPModAS}$ and let $(X,a)$ be a $\nabla^I$-category. By Example \ref{exmp:prodq}.(1), $(X,a_i)$ is a $\nabla$-category for all $i\in I$, where $a_i$ is the $i$th-coordinate function of $a.$ Hence $\big\{(X,w_{a_i}):i\in I\big\}$ is a family of quasi-pseudometric modular spaces (see Theorem \ref{thm:iso}).
	By assumption $(X,F\circ w_\Delta)$ is a quasi-pseudometric modular space so $(X,a_{F\circ w_\Delta})$ is a $\nabla$-category  by Theorem \ref{thm:iso}. 
	Notice that, given $x,y\in X$ and $t>0$
	\begin{align*}
		(a_{F\circ w_\Delta}(x,y))(t)&=(F\circ w_\Delta)(t,x,y)=F\big((w_{a_i}(t,x,y))_{i\in I}\big)\\&=F_\nabla\Big(\big(w_{a_i}(\cdot,x,y)\big)_{i\in I}\Big)(t)=F_\nabla\Big(\big(a_i(x,y)\big)_{i\in I}\Big)(t)\\
		&=\big((F_\nabla\circ a)(x,y)\big)(t)
	\end{align*}
	that is,
	$$a_{F\circ w_\Delta}=F_\nabla\circ a,$$
	that proves the implication.
	
	(3) $\Rightarrow$ (1) Let $\big\{(X_i,w_i):i\in I\big\}$ be a family of quasi-pseudometric modular spaces.  By Theorem \ref{thm:iso}, $\big\{(X_i,a_{w_i}):i\in I\big\}$ is a family of $\nabla$-categories so $\left(\prod_{i\in I} X_i,a_\Pi\right)$ is $\nabla^I$-category where, for all $i\in I$,
	$$\big(a_\Pi(x,y)\big)_i(t)=\big(a_{w_i}(x_i,y_i)\big)(t)=w_i(t,x_i,y_i)$$
	for every $x,y\in \prod_{i\in I} X_i$, $t>0$ (see Example \ref{exmp:prodq}.(2)). By assumption, $\Big(\prod_{i\in I} X_i, F_\nabla\circ a_\Pi\Big)$ is a $\nabla$-category so, by Theorem \ref{thm:iso}, $\Big(\prod_{i\in I} X_i, w_{F_\nabla\circ a_\Pi}\Big)$ is a quasi-pseudometric modular space. Moreover, given $x,y\in\prod_{i\in I} X_i$ and $t>0$
	\begin{align*}
		w_{F_\nabla\circ a_\Pi}(t,x,y)&=\Big(\big(F_\nabla\circ a_\Pi\big)(x,y)\Big)(t)=\Big(F_\nabla\big((w_i(\cdot,x_i,y_i))_{i\in I}\big)\Big)(t)\\&=F\Big((w_i(t,x_i,y_i)_{i\in I}\Big)
		=(F\circ w_\Pi)(t,x,y)
	\end{align*}
	so $\left(\prod_{i\in I}X_i,w_{F_\nabla\circ a_w}\right)=\left(\prod_{i\in I}X_i,F\circ w_\Pi\right)$ is a quasi-pseudometric modular space. Hence $F\in\mathsf{QPModAP}.$
\end{proof}

\begin{remark}
	We observe that given a function $F: \left([0,+\infty]^{(0,+\infty)}\right)^I\to [0,+\infty]^{(0,+\infty)}$ satisfying that $(X,F\circ a)$ is a $\nabla$-category for every $\nabla^I$-category $(X,a),$ then
	$$F(\nabla^I)\subseteq \nabla.$$
	Let us show this. Let $(f_i)_{i\in I}\in\nabla^I.$ Consider a set $X=\{x_1,x_2\}$ with two different points and for every $i\in I$, let $a_{f_i}:X\times X\to\nabla$ be defined as
	$$\big(a_{f_i}(x,y)\big)(t)=\begin{cases}
		f_i(t)&\text{ if }x\neq y,\\
		0&\text{ if }x=y,
	\end{cases}$$
	for every $x,y\in X,$ $t>0.$  Then $(X,a_f)$ is $\nabla^I$-category where
	$$\big(a_f(x,y)\big)_i(t)=\big(a_{f_i}(x,y)\big)(t)$$
	for every $x,y\in X,$ $t>0,$ $i\in I.$ By assumption $(X,F\circ a_f)$ is a $\nabla$-category so 
	$$(F\circ a_f)(x_1,x_2)=F\big((a_f(x_1,x_2))_{i\in I}\big)=F\big((f_i)_{i\in I}\big)\in\nabla$$
	which proves the assertion.
\end{remark}

Based on the previous proposition and remark, we can derive a corollary that shows quasi-pseudometric modular aggregation functions are included in the family of preserving functions between $\nabla^I$ and $\nabla$. This connects the existing theory of quasi-pseudometric modular aggregation functions in the literature with our approach to the problem using preserving functions.

\begin{corollary}\label{cor:modaggispres}
	Let $F:[0,+\infty]^I\to [0,+\infty]$ be a function. Define $F_\nabla: \left([0,+\infty]^{(0,+\infty)}\right)^I\to [0,+\infty]^{(0,+\infty)}$ as
	$$F_\nabla\left((f_i)_{i\in I}\right)(t)=F((f_i(t))_{i\in I})$$
	for every $(f_i)_{i\in I}\in \left([0,+\infty]^{(0,+\infty)}\right)^I$ and $t>0.$ 
	Then the following statements are equivalent:
	\begin{enumerate}
		\item $F\in \mathsf{QPModAP};$
		\item $F\in \mathsf{QPModAS};$
		\item $F_\nabla\in\mathscr{P}(\nabla^I,\nabla).$
	\end{enumerate}
	Consequently, the map
	$$\nu: \mathsf{QPModAP} \to \mathscr{P}(\nabla^{I},\nabla)$$
	given by
	$$\nu(F)=F_{\nabla}$$
	is well-defined.
\end{corollary}

In \cite[Theorem 8]{BFVale24}, Bibiloni-Femenias and Valero characterized the quasi-pseudometric modular aggregation functions on sets as those isotone and subadditive functions $F:[0,+\infty]^I\to [0,+\infty]$ such that $F(0_I)=0.$ We will show that this result can be derived from the results presented in \cite{FBRL}. To achieve this, we will prove that a function $F:\mathsf{P}_+^I\to \mathsf{P}_+$ is preserving if and only if $F_{\nabla}:\nabla^{I}\to \nabla$ is preserving. This will be a consequence of the following result.

\begin{proposition}\label{prop:modlaxmor}
	Let $F:\mathsf{P}_+^I\to \mathsf{P}_+$. Then $F$ is a lax morphism of quantales if and only if $F_{\nabla}:\nabla^I\to \nabla$ is a lax morphism of quantales.
\end{proposition}

\begin{proof}
	Suppose that $F$ is a lax morphism of quantales.
	Since $F$ is a lax morphism then it is isotone. Hence $F_\nabla(\nabla^I)\subseteq \nabla$ so $F_\nabla$ is well-defined.
	
	Next, let us check that $F_{\nabla}$ is isotone. Consider $(f_{i})_{i\in I}, (g_{i})_{i\in I} \in \nabla^{I}.$ Since $F$ is subadditive and isotone (Example \ref{ex:lax_morphism_P}) then
	\begin{align*}
		\Big(F_{\nabla}\big((f_{i})_{i\in I}\big) \oplus F_{\nabla}\big((g_{i})_{i\in I}\big)\Big)(t)&=\bigvee_{u+v=t}^{\leq^{op}}F\big((f_{i}(u))_{i\in I}\big)+F\big((g_{i}(v))_{i\in I}\big)\leq^{op}\\
		&\leq^{op}\bigvee_{u+v=t}^{\leq^{op}}F\big((f_{i}(u)+g_{i}(v))_{i\in I}\big)\\&\leq^{op} F\left(\left(\bigvee_{u+v=t}^{\leq^{op}}\big(f_{i}(u)+g_{i}(v)\big)\right)_{i\in I}\right)\\&=\Big(F_{\nabla}\big((f_{i}\oplus g_{i})_{i\in I}\big)\Big)(t).
	\end{align*}
	
	Hence, $F_\nabla$ is subadditive.
	
	Moreover, suppose that $f_{i} \leq^{op} g_{i}$ for all $i\in I$. Using that $F$ is isotone we have:
	$$F_{\nabla}\big((f_{i})_{i\in I}\big)(t)=F\big((f_{i}(t)\big)_{i\in I})\leq^{op} F\big((g_{i}(t))_{i\in I}\big)=F_{\nabla}\big((g_{i})_{i\in I}\big)(t).$$
	Finally, we have that, for every $t>0$
	$$(F_{\nabla}(0_{\nabla^I}))(t)=F(0_I)=0$$
	so $F_{\nabla}(0_{\nabla^I})=0_\nabla.$
	Therefore, $F_\nabla$ is a lax morphism of quantales.
	
	Conversely, suppose that $F_\nabla$ is a lax morphism of quantales. 
	We first notice that, given $t>0$
	$$F(0_I)=F(0_{\nabla^I} (t))=\big(F_\nabla(0_{\nabla^I})\big)(t)=0_\nabla (t)=0.$$
	Now we prove that $F$ is subadditive.\\
	Let $x=(x_{i})_{i \in I},y=(y_{i})_{i \in I}\in [0,+\infty]^I$. For each $i\in I$, let us consider $f_i,g_i\in \nabla$ defined as 
	\begin{align*}
		f_{i}(t)&=x_{i},\\
		g_{i}(t)&=y_{i},
	\end{align*}
	for every $t\in (0,+\infty).$ Let $f:=(f_i)_{i\in I}, g:=(g_i)_{i\in I}\in\nabla^I.$ Then, for any $t>0,$ using the subadditivity of $F_\nabla$, we deduce that, for any $t>0$, 
	\begin{align*} 			
		F(x+y)&=F\big((f_{i}(t)+g_{i}(t))_{i \in I}\big)=\big(F_{\nabla}((f_{i}+g_{i})_{i \in I})\big)(t)=\big(F_\nabla(f\oplus g)\big)(t)\\
		&\leq \Big( \big(F_\nabla(f)\big)\oplus \big(F_\nabla(g)\big) \Big)(t)= \big(F_{\nabla}(f)\big)(t)+\big(F_{\nabla}(g)\big)(t)\\&=F\big((f_i(t))_{i\in I}\big)+F\big((g_i(t))_{i\in I}\big)=F(x)+F(y).
	\end{align*}
	so $F$ is subadditive. 
	
	Furthermore, if $x\leq y$ then $g\leq^{\text{op}} f$ and since $F_\nabla$ is isotone we have
	$$F(x)=F\big((f_i(t))_{i\in I}\big)=\big(F_\nabla(f)\big)(t)\leq \big(F_\nabla(g)\big)(t)=F\big((g_i(t))_{i\in I}\big)=F(y).$$
	Consequently, $F$ is isotone and a lax morphism of quantales.	
\end{proof}

Based on our previous results, we can achieve one of the main goals of the paper: the characterization of the quasi-pseudometric modular aggregation functions on products  or on sets. The following theorem characterizes them in various ways.

\begin{theorem}[\mbox{compare with  \cite[Theorem 8]{BFVale24}}]\label{thm:FpresifFnablapres}
	Let $F:\mathsf{P}_+^I\to \mathsf{P}_+$. Then the following assertions are equivalent:
	\begin{enumerate}
		\item $F \in \mathsf{QPModAP}$;
		\item $F \in \mathsf{QPModAS}$;
		\item $F_{\nabla}\in\mathscr{P}(\nabla^I,\nabla)$;
		\item $F\in \mathscr{P}(\mathsf{P}_+^I,\mathsf{P}_+)$;
		\item $F(0_I)=0$, $F$ is isotone and $F$ is subadditive;
		\item $F(0_I)=0$ and $F$ preserves asymmetric triangle triplets;
		\item $F(0_I)=0$, $F$ is isotone and $F$ preserverse triangle triplets.
	\end{enumerate}
	
\end{theorem}

\begin{proof}
	This follows from Corollary \ref{cor:modaggispres}, Proposition \ref{prop:modlaxmor} and Theorem \ref{thm:plmq}.
\end{proof}

\begin{remark}
	Notice that the equivalence between (2), (5), (6), and (7) was proved in \cite[Theorem 8]{BFVale24}. Nevertheless, we have been able to obtain this equivalence based on the general theory regarding preserving functions. This approach also enables us to enhance the result by providing the same characterization for quasi-pseudometric modular aggregation functions on products.
\end{remark}

\begin{example}\label{ex:af}
	Let $(k_n)_{n\in\mathbb{N}}$ be a sequence of positive real numbers. Consider the following functions:
	\begin{enumerate}
		\item $F:[0,+\infty]^\mathbb{N}\to[0,+\infty]$ given by
			$$F(x)=\sum_{n\in\mathbb{N}}k_n\cdot x_n$$
			for all $x=(x_n)_{n\in\mathbb{N}}$;
		\item $G:[0,+\infty]^\mathbb{N}\to[0,+\infty]$ given by
			$$G(x)=\sup\{k_n\cdot x_n:n\in\mathbb{N}\}$$
			for all $x=(x_n)_{n\in\mathbb{N}}.$
	\end{enumerate}
	Then $F,G$ are quasi-pseudometric modular aggregation functions on products and on sets. This is a direct consequence of the previous theorem since $F(0_\mathbb{N})=G(0_\mathbb{N})=0$, and $F,G$ are clearly isotone and subadditive.
\end{example}

Based on our previous discussion, we can identify two distinct types of functions for aggregating quasi-pseudometric modulars. One type consists of functions of the form \( F: [0, +\infty]^I \to [0, +\infty] \), while the other involves preserving functions between the quantales \( \nabla^I \) and \( \nabla \). This leads us to the following definition.

\begin{definition}
	Let $F:\nabla^{I} \to \nabla$. We say that $F$ is a
	\begin{itemize}
		\item \textbf{(quasi-)(pseudo)metric modular aggregation function on pro\-ducts} if for every family $\big\{ (X_{i}, w_{i}) \text{ : } i \in I\big\}$ of (quasi-)(pseudo)metric modular spaces then $(\Pi_{i \in I}X_{i}, w_\Pi^F)$ is a (quasi-)(pseudo)metric modu\-lar space, where
		$$w_\Pi^F(t,x,y)=\Big(F\big((w_{i}(\cdot,x_{i},y_{i}))_{i}\big)\Big)(t)$$
		for every $x,y\in \prod_{i\in I} X_i,t>0.$
		\item \textbf{(quasi-)(pseudo)metric modular aggregation function on sets} if for every collection $\big\{ w_{i}:i\in I \big\}$ of (quasi-)(pseudo)metric modulars over a set $X$, we have that $(X,w_\Delta^F)$ is a  (quasi-)(pseudo)metric modular space, where
		$$w_{\Delta}^F(t,x,y)=\Big(F\big((w_{i}(\cdot,x,y))_{i}\big)\Big)(t)$$
		for every $x,y\in X,t>0.$
	\end{itemize}
\end{definition}

The next result shows that the preserving functions between the quantales $\nabla^I$ and $\nabla$ are precisely the quasi-pseudometric modular aggregation functions as defined above.

\begin{proposition}\label{prop:preservingnabla}
	Let $F:\nabla^{I} \to \nabla$. Then the following statements are equivalent:
	\begin{enumerate}
		\item $F\in\mathscr{P}(\nabla^I,\nabla);$
		\item $F$ is a quasi-pseudometric modular aggregration function on pro\-ducts;
		\item $F$ is a quasi-pseudometric modular aggregration function on sets.
	\end{enumerate}
\end{proposition}

\begin{proof}
	(1) $\Rightarrow$ (2) Let $\big\{ (X_{i}, w_{i}) \text{ : } i \in I\big\}$ be a family of quasi-pseudometric modular spaces. Then $(\prod_{i\in I} X_i, a_\Pi^w)$ is a $\nabla^I$-category where $(a_\Pi^w (x,y))_i(t)=w_i(t,x_i,y_i)$ for all $x,y\in\prod_{i\in I}X_i, t>0.$ By assumption, $(\prod_{i\in I} X_i, F\circ a_\Pi^w)$ is a $\nabla$-category so $(\prod_{i\in I} X_i, w_{F\circ a_\Pi^w})$ is a quasi-pseudometric modular space (see Theorem \ref{thm:iso}). Since
	\begin{align*}
		w_{F\circ a_\Pi^w}(t,x,y)&=(F\circ a_\Pi^w(x,y))(t)=\Big(F\big((w_{i}(\cdot,x_{i},y_{i}))_{i\in I}\big)\Big)(t)=w_\Pi^F(t,x,y)
	\end{align*}
	for all $x,y\in\prod_{i\in I}X_i,t>0$ then $(\Pi_{i \in I}X_{i}, w_\Pi^F)$ is a quasi-pseudometric modular space. Hence, $F$ is a quasi-pseudometric modular aggregation function on products.
	
	(2) $\Rightarrow$ (3) This is straightforward.
	
	(3) $\Rightarrow$ (1) Let $(X,a)$ be a $\nabla^I$-category. Then $\big\{(X,w_{a_i}):i\in I\big\}$ is a family of quasi-pseudometric modular spaces (see Example \ref{exmp:prodq}.(1) and Theorem \ref{thm:iso}). By hypothesis, $(X,w_\Delta^F)$ is a quasi-pseudometric mo\-dular space, so $(X,a_{w_\Delta^F})$ is a $\nabla$-category. Moreover,
	\begin{align*}
		(a_{w_\Delta^F}(x,y))(t)&=w_\Delta^F(t,x,y)=\Big(F\big((w_{a_i}(\cdot,x,y))_{i\in I}\big)\Big)(t)=\Big(F\big((a_i(x,y))_{i\in I}\big)\Big)(t)\\
		&=\Big(F\big(a(x,y)\big)\Big)(t)=((F\circ a)(x,y))(t)
	\end{align*}
	for all $x,y\in X, t>0.$ Consequently, $(X,F\circ a)$ is a $\nabla$-category.
	
\end{proof}

From the previous proposition and Theorem \ref{thm:FpresifFnablapres}, we can deduce a result that connects the families $\mathscr{P}(\mathsf{P}_+^I,\mathsf{P}_+)$ and $\mathscr{P}(\nabla^I,\nabla).$

\begin{proposition}\label{prop:presqpmm}
	Let $F:\mathsf{P}_+^I\to \mathsf{P}_+$. The following statements are equivalent.
	\begin{enumerate}
		\item $F$ is a quasi-pseudometric modular aggregation function on products;
		\item $F$ is a quasi-pseudometric modular aggregation function on sets;
		\item $F\in\mathscr{P}(\mathsf{P}_+^I,\mathsf{P}_+)$;
		\item $F_{\nabla}$ is a quasi-pseudometric modular aggregation function on pro\-ducts;
		\item $F_{\nabla}$ is a quasi-pseudometric modular aggregation function on sets;
		\item $F_\nabla\in\mathscr{P}(\nabla^I,\nabla)$.
	\end{enumerate}
	
\end{proposition}

To summarize, in this section we have presented one of the main contributions of the paper: the characterization of quasi-pseudometric modular aggregation functions on products, which we have demonstrated to be equivalent to aggregation on sets. Furthermore, we have established a relationship between the classical method for aggregating quasi-pseudometric modulars, as developed in \cite{BFMinyaVale23,BFVale24}, and preserving functions between certain quantales.
\section{Aggregation of pseudometric modulars}\label{sec:apmod}

Bibiloni-Femenias, Mi\~nana and Valero characterized in \cite[Theorem 5]{BFMinyaVale23} the pseudometric modular aggregation functions on sets.  In this section, we address the problem of deriving their result using our general framework, as well as characterizing pseudometric modular aggregation functions on products. 

The following result highlights two key points. First, the family of pseudometric modular aggregation functions on products is equal to the family of pseudometric modular aggregation functions on sets. Second, symmetrically preserving functions (see Definition \ref{def:preserving}) are suitable for interpreting pseudometric modular aggregation functions within our theoretical framework.

\begin{proposition}\label{prop:pmod_syP}
	Let $F:\mathsf{P}_+^I\to \mathsf{P}_+$. The following statements are equivalent.
	\begin{enumerate}
		\item $F\in \mathsf{PModAP}$;
		\item $F\in \mathsf{PModAS}$;
		\item $F_\nabla\in\mathsf{sy}\mathscr{P}(\nabla^I,\nabla).$
	\end{enumerate}
\end{proposition}

\begin{proof}
	The proof is similar to that of Corollary \ref{cor:modaggispres} taking into account Remark \ref{rem:symsep}.
\end{proof}

The next result will allow us to prove the characterization of pseudometric modular aggregation functions.

\begin{proposition}\label{prop:symsiso}
	Let $F:\mathsf{P}_+^I\to \mathsf{P}_+.$ Then $F_\nabla\in\mathsf{sy}\mathscr{P}(\nabla^I,\nabla)$ if and only if $F\in\mathsf{sy}\mathscr{P}(\mathsf{P}_+^I,\mathsf{P}_+)$ and $F$ is isotone.
\end{proposition}

\begin{proof}
	Suppose that $F_\nabla$ is symmetrically preserving. Let $(X,a)$ be a symmetric  $\mathsf{P}_+^I$-category. Hence $\big\{(X,a_i):i\in I\big\}$ is a family of symmetric $\mathsf{P}_+$-categories (Example \ref{exmp:prodq}), that is, a family of extended pseudometric spaces (Example \ref{ex:vcat}). Therefore  $\big\{(X,w_{a_i}):i\in I\big\}$ is a family of pseudometric modular spaces (Example \ref{exmp:wgd}) so $\big\{(X,a_{w_{a_i}}):i\in I\big\}$ is a family of symmetric $\nabla$-categories (Theorem \ref{thm:iso}). By Example \ref{exmp:prodq}.(3), $(X,a_\Delta)$ is a symmetric $\nabla^I$-category so, by assumption, $(X,F\circ a_\Delta)$ is a symmetric $\nabla$-category. Observe that
	\begin{align*}
		\Big(\big(F\circ a_\Delta\big)(x,y)\Big)(t)&=F\Big(\big(\big(a_{w_{a_i}}(x,y)\big)(t)\big)_{i\in I}\Big)=F\Big(\big(w_{a_i}(t,x,y)\big)_{i\in I}\Big)\\
		&=F\Big(\big(a_i(x,y)\big)_{i\in I}\Big)=(F\circ a)(x,y)
	\end{align*}
	for all $t>0.$ Hence, $(X,F\circ a)$ is a symmetric $\mathsf{P}_+$-category so $F\in\mathsf{sy}\mathscr{P}(\mathsf{P}_+^I,\mathsf{P}_+).$
	
	Next, we will prove that $F$ is isotone. Notice that, by Theorem \ref{thm:symp}, $F_\nabla(0_{\nabla^I})=0_\nabla$ so $F(0_I)=0.$ 
	
	Let $s,t\in [0,+\infty]^I$ such that $s\leq t.$ Consider a set $X=\{0,1\}$ with two different elements and, for each $i\in I$ define $a_i:X\times X\to \nabla$ as
	$$\big(a_i(x,y)\big)(t)=\begin{cases}
		0&\text{ if }x=y, t>0,\\
		t_i&\text{ if }x\neq y, 0<t<1,\\
		s_i&\text{ if }x\neq y, 1\leq t,
	\end{cases}$$
	for all $x,y\in X$, $t>0.$ It is straightforward to check that $(X,a_i)$ is a symmetric $\nabla$-category for all $i\in I.$ Then $(X,a_\Delta)$ is a symmetric $\nabla^I$-category (Example \ref{exmp:prodq}.(3)) so $(X,F_\nabla\circ a_\Delta)$ is a symmetric $\nabla$-category. Then
	\begin{align*}
		F(s)&=F\Big(\big((a_i(0,1))(1)\big)_{i\in I}\Big)=F\Big(\big(a_\Delta(0,1)\big)(1)\Big)=\Big(\big(F_\nabla\circ a_\Delta\big)(0,1)\Big)(1)\\
		&\leq \left(\Big(\big(F_\nabla\circ a_\Delta\big)(0,1)\Big)\oplus \Big(\big(F_\nabla\circ a_\Delta\big)(1,1)\Big)\right)(1)\\&=\bigwedge_{0\leq s\leq 1} \Big(\big(F_\nabla\circ a_\Delta\big)(0,1)\Big)(s)+ \Big(\big(F_\nabla\circ a_\Delta\big)(1,1)\Big)(1-s) \\
		&\leq \Big(\big(F_\nabla\circ a_\Delta\big)(0,1)\Big)(\tfrac{1}{2})+ \Big(\big(F_\nabla\circ a_\Delta\big)(1,1)\Big)(\tfrac{1}{2})\\&=
		F\Big(\big((a_i(0,1))(\tfrac{1}{2})\big)_{i\in I}\Big)+F\Big(\big((a_i(0,0))(\tfrac{1}{2})\big)_{i\in I}\Big)=F(t)+F(0_I)=F(t)
	\end{align*}
	so $F$ is isotone.
	
	Conversely, suppose that $F\in\mathsf{sy}\mathscr{P}(\mathsf{P}_+^I,\mathsf{P}_+)$ and $F$ is isotone. By Theorem \ref{thm:symp} $F$ preserves triangle triplets and $F(0_I)=0$. Given $x,y\in [0,+\infty]^I$ then $(x+y,x,y)$ is a triangle triplet on $[0,+\infty]^I$ so $(F(x+y),F(x),F(y))$ is a triangle triplet on $[0,+\infty]$. Hence $F(x+y)\leq F(x)+F(y)$ so $F$ is subadditive. Therefore, $F$ is a lax morphism of quantales which implies, by Proposition \ref{prop:modlaxmor} and Theorem \ref{thm:plmq} that $F_\nabla \in \mathscr{P}(\nabla^I,\nabla) \subseteq \mathsf{sy}\mathscr{P}(\nabla^I,\nabla).$
\end{proof}

The following result shows that there is no difference between the quasi-pseudometric modular aggregation functions (on products or on sets) and the pseudometric modular aggregation functions (on products or on sets).

\begin{theorem}[\mbox{compare with \cite[Theorem 8]{BFVale24}}]\label{thm:pmod_long}
	Let $F:\mathsf{P}^{I}_+ \to \mathsf{P}$. The following statements are equivalent.
	\begin{enumerate}
		\item $F\in \mathscr{P}(\mathsf{P}_+^I,\mathsf{P}_+)$;
		\item $F\in \mathsf{sy}\mathscr{P}(\mathsf{P}_+^I,\mathsf{P}_+)$ and $F$ is isotone;
		\item $F_{\nabla}\in\mathscr{P}(\nabla^I,\nabla)$;
		\item $F_{\nabla}\in\mathsf{sy}\mathscr{P}(\nabla^I,\nabla)$; 
		\item $F \in \mathsf{QPModAP}$;
		\item $F \in \mathsf{QPModAS}$;
		\item $F_{\nabla}$ is a quasi-pseudometric modular aggregation function on products;
		\item $F_{\nabla}$ is a quasi-pseudometric modular aggregation function on sets;
		\item $F \in \mathsf{PModAP}$;
		\item $F \in \mathsf{PModAS}$;
		\item $F(0_I)=0$, $F$ is isotone and $F$ is subadditive;
		\item $F(0_I)=0$ and $F$ preserves asymmetric triangle triplets;
		\item $F(0_I)=0$, $F$ is isotone and $F$ preserves triangle triplets.
	\end{enumerate}
\end{theorem}

\begin{proof}
	The equivalence between (1) and (2) follows from Theorems \ref{thm:plmq}, \ref{thm:symp} and Proposition \ref{prop:asym_symiso}. Moreover, (1) is equivalent to (3), (5), (6), (7) and (8) by Proposition \ref{prop:presqpmm}. All these statements are equivalent to (11), (12) and (13) by Theorem \ref{thm:FpresifFnablapres}.
	
	(2) is equivalent to (4) by Proposition \ref{prop:symsiso}. (4) is equivalent to (9) and (10) by Proposition \ref{prop:pmod_syP}.

\end{proof}

\begin{remark}
	Notice that the equivalence between (6), (10), (11), (12) and (13) was first proved in \cite[Theorem 8]{BFVale24}.
\end{remark}

\begin{remark}
	Observe that by \cite[Theorem 3.15]{FBRL}, the statements of Theorem \ref{thm:pmod_long} are equivalent to $F$ is an extended quasi-pseudometric aggregation function on products or on sets.
	
	Moreover, although in the context of metric modulars, the concepts of pseudometric modular aggregation functions (on products or on sets) and quasi-pseudometric modular aggregation functions (on products or on sets) are equivalent, this equivalence does not hold for metrics. There exist pseudometric aggregation functions on sets that are not quasi-pseudometric aggregation functions on sets (see \cite{MayorVale10}).
\end{remark}

In this section, we have continued the work developed in the previous one, characterizing aggregation functions of pseudometric modular, showing that they coincide with their quasi-pseudometric modular counterparts. Using the categorical framework developed earlier, we have demonstrated that symmetrically preserving functions between certain quantales provide an alternative interpretation of pseudometric modular aggregation functions.

\section{Aggregation of quasi-metric modulars}\label{sec:aqmod}

We now address the problem of characterizing the quasi-metric modular aggregation functions, both on products or on sets. This question was solved in \cite[Theorem 9]{BFVale24} when the aggregation is performed on sets. 

Next, we will demonstrate that the separately preserving functions between $\mathsf{P}_+^I$ and $\mathsf{P}_+$ are precisely the quasi-metric modular aggregation functions on products. We also provide an internal characterization of these functions.

\begin{theorem}\label{thm:sepp_qmmodp}
	Let $F:\mathsf{P}_+^I\to \mathsf{P}_+$. The following statements are equivalent:
	\begin{enumerate}
		\item $F\in \mathsf{sep}\mathscr{P}(\mathsf{P}_+^I,\mathsf{P}_+);$
		\item $F\in\mathsf{QModAP}$; 
		\item $F$ is isotone, subadditive and $F^{-1}(0)=\{0_I\}.$
	\end{enumerate}
\end{theorem}

\begin{proof}
	(1) $\Rightarrow$ (2) Let us suppose that $F$ is separately preserving. Let $\big\lbrace (X_{i}, w_{i}):i\in I \big\rbrace$ be a family of quasi-metric modular spaces. By Theorems \ref{thm:plmq} and \ref{thm:seplmq},  $F$ is preserving so it is a quasi-pseudometric modular aggregation function on products by Theorem \ref{thm:FpresifFnablapres}. Hence $(\prod_{i\in I} X_i, F\circ w_\Pi)$ is a quasi-pseudometric modular space. Consequently, it only remains to show the property (M3) from Definition \ref{def:mqpm}.  Let $x,y\in\prod_{i\in I}{X_i}$ such that $F\circ w_\Pi(t,x,y)=F\circ w_\Pi(t,y,x)=0$ for every $t>0$, that is, $F\big((w_{i}(t,x_{i},y_{i}))_{i\in I}\big)=F\big((w_{i}(t,y_{i},x_{i}))_{i\in I}\big)=0$. Since $F$ is separately preserving then $F^{-1}(0)=0_I$ by Theorem \ref{thm:seplmq}, so $w_{i}(t,x_{i},y_{i})=0$ for every $t >0$ and $i \in I.$ Since $w_{i}$ is a quasi-metric modular for every $i \in I$,  we conclude $x_{i}=y_{i}$ for every $i \in I$, as desired.
	
	(2) $\Rightarrow$ (3) If $F$ is a quasi-metric modular aggregation function on pro\-ducts then it is a quasi-metric modular aggregation function on sets. By \cite[Theorem 9]{BFVale24}, $F$ is isotone, subadditive, and $F(0_I)=0.$
	
	Let us show that $F^{-1}(0)=\{0_I\}.$ By way of contradiction, suppose that there exists $a\in F^{-1}(0)$ such that $a\neq 0_I$. Consider $\big([0,+\infty),w\big)$ the metric modular space of Example \ref{exmp:mods}. Then $\{\big([0,+\infty),w\big):i\in I\}$ is a family of quasi-metric modular spaces so, by assumption, $\big([0,+\infty]^I,F\circ w_\Pi\big)$ is a quasi-metric modular space. Nevertheless, $a\neq 0_I$ but
	$$(F\circ w_\Pi)(t,a,0_I)=F((w_i(t,a_i,0))_{i\in i})=F((a_i)_{i\in I})=0=(F\circ w_\Pi)(t,0_I,a)$$
	for all $t>0.$ This contradicts axiom (M3) of a quasi-metric modular.

	(3) $\Rightarrow$ (1) This is a consequence of Theorem $\ref{thm:seplmq}$ .
\end{proof}

\begin{remark}
	According to \cite[Theorem 3.18]{FBRL}, the statements of the previous theorem are equivalent to $F$ is an extended quasi-metric aggregation function on products.
\end{remark}

Notice that all functions in Example \ref{ex:af} are quasi-metric modular aggregation functions on products, since they satisfy statement 3 in Theorem \ref{thm:sepp_qmmodp}. The constant zero function between $\mathsf{P}_+^I$ and $\mathsf{P}_+$ is an example of a quasi-pseudometric modular aggregation function on products which is not a quasi-metric modular aggregation function on products.

As previously mentioned, Bibiloni-Femenias and Valero characterized in \cite{BFVale24} the quasi-metric modular aggregation functions on sets through the following result:

\begin{theorem}[\mbox{\cite[Theorem 9]{BFVale24}}]\label{thm:BFVAle24}
	Let $n \in \mathbb{N}$, and let $F:[0,+\infty]^n \rightarrow[0,+\infty]$ be a function. The following statements are equivalent to each other:
	\begin{enumerate}
		\item $F$ is a quasi-metric modular aggregation function on sets;
		\item $F$ is a metric modular aggregation function on sets;
		\item $F\left(0_n\right)=0$ and $F$ is isotone and subadditive. Moreover, if $a \in[0,+\infty]^n$ and $F(a)=0$, then $a_i=0$ for some $i=1, \ldots, n$;
		\item $F\left(0_n\right)=0$, and, in addition, $F(c) \leq F(a)+F(b)$ for all $a, b, c \in[0,+\infty]^n$, with $c \leq a+b$. Moreover, if $a \in[0,+\infty]^n$ and $F(a)=0$, then $a_1=0$ for some $i=1, \ldots, n$;
		\item $F\left(0_n\right)=0, F$ is isotone and preserves triangular triplets. Moreover, if $a \in[0,+\infty]^n$ and $F(a)=0$, then $a_1=0$ for some $i=1, \ldots, n$.
	\end{enumerate}
\end{theorem}

Notice that from the above Theorem and Theorem \ref{thm:sepp_qmmodp} we deduce that every quasi-metric modular aggregation function on products is a quasi-metric modular aggregation function on sets, but the converse does not hold as the next trivial example shows.

\begin{example}\label{ex:sets_not_products}
	Given $j\in I,$ let us consider $P_{j}:[0,+\infty]^{I} \to [0,+\infty]$ be the $j$th projection, that is, 
	$$P_{j}(x)=x_{j}$$
	for all $x=(x_{i})_{i \in I}\in [0,+\infty]^{I}.$ Let us check that $P_{j}$ is a quasi-metric modular aggregation function on sets. Indeed, consider a collection $\lbrace w_{i} \text{ : } i \in I \rbrace$ of quasi-metric modulars on a fixed nonempty set $X.$  For every $x, y \in X$ and $t>0,$
		$$P_{j}\Big(\big(w_{i}(t,x,y)\big)_{i \in I}\Big)=w_{j}(t,x,y),$$ that is, $P_j\circ w_\Delta=w_j$, which is a quasi-metric modular on $X.$ Thus $P_{j}$ is a quasi-metric modular aggregation function on sets. 
	
	However, $P_j$ is not a quasi-metric modular  aggregation function on products. To check this, let $(X,w)$ be an arbitrary quasi-metric modular space such that $X$ has at least two different points $x,y$. Then $(X^I, F\circ w_\Pi)$ is not a quasi-metric modular space since if we consider $\chi,\mu\in X^I$ such that $\chi_i=x$ for all $i\in I$ and
	$$\mu_i=\begin{cases}
		y&\text{ if }i\neq j,\\
		x&\text{ if }i=j,
	\end{cases}$$
	then 
	$$(P_j\circ w_\Pi)(t,\chi,\mu)=P_j((w(t,\chi_i,\mu_i))_{i\in I})=w(t,\chi_j,\mu_j)=w(t,x,x)=0$$
	but $\chi\neq\mu,$ so $P_j\circ w_\Pi$ is not a quasi-metric modular on $X^I.$
	
\end{example}

We obtain the following characterization of functions that aggregate metric modulars on products based on the previous results. We must note that we prove these are equivalent to the quasi-metric modular aggregation functions on products.

\begin{theorem}
	Let $F:\mathsf{P}_+^I\to \mathsf{P}_+$. The following statements are equivalent:
	\begin{enumerate}
		\item $F\in \mathsf{sep}\mathscr{P}(\mathsf{P}_+^I,\mathsf{P}_+)$;
		\item $F\in\mathscr{P}(\mathsf{P}_+^I,\mathsf{P}_+)$ and $F^{-1}(0)=\{0_I\}$;
		\item $F\in\mathsf{sym}\mathscr{P}(\mathsf{P}_+^I,\mathsf{P}_+),$ $F$ is isotone and $F^{-1}(0)=\{0_I\}$;
		\item $F\in\mathsf{ModAP}$; 
		\item $F$ is isotone, subadditive and $F^{-1}(0)=\{0_I\};$
		\item $F \in \mathsf{QModAP}.$
	\end{enumerate}
\end{theorem}

\begin{proof}
	(1) $\Leftrightarrow$ (2) $\Leftrightarrow$ (3) These equivalences follow from Theorems \ref{thm:sepp_qmmodp} and \ref{thm:pmod_long}.
	
	(3) $\Rightarrow$ (4) Since $F$ is symmetrically preserving and isotone then, by Theorem \ref{thm:pmod_long}, it is a pseudometric modular aggregation function on products. Moreover, since $F$ is separately preserving, then $F$ is a quasi-metric modular aggregation function on products by Theorem \ref{thm:sepp_qmmodp}. This obviously implies that $F$ is a modular metric aggregation function on products.
	
	(4) $\Rightarrow$ (5) This is similar to (2) $\Rightarrow$ (3) of Theorem \ref{thm:sepp_qmmodp}. 
	
	(5) $\Rightarrow$ (6) $\Rightarrow$ (1) This is a consequence of Theorem \ref{thm:sepp_qmmodp}.
\end{proof}

In conclusion, this section has focused on characterizing functions that aggregate quasi-metric modulars. However, unlike previous cases, the results differ depending on whether the aggregation occurs on products or on sets. Additionally, we have demonstrated that quasi-metric modular aggregation functions on products coincide with metric modular aggregation functions on products. This fulfills our objective of characterizing the functions aggregating on products these mathematical structures.

\section{Conclusions and future work}

In this paper, we have investigated functions that aggregate quasi-pseudo\-metric modulars. While this topic has been recently studied in the literature \cite{BFMinyaVale23,BFVale24}, it has only focused on aggregation on sets. Our goal is to fill a gap in this research by examining aggregation on products. Our approach builds on recent advancements in the theory of aggregation functions over quantales \cite{FBRL}, emphasizing the relevance and applicability of this theory.

A significant contribution of this work is the characterization of aggregation functions for quasi-pseudometric modulars on products. We have provided an internal characterization of these functions. Moreover,  we have proved that the aggregation functions for quasi-pseudometric modulars on products are exactly the same as those for quasi-pseudometric modulars on sets, revealing an unexpected structural parallel between these two approaches. This similarity also occurs for pseudometric modulars. However, it does not hold for quasi-metric modulars or metric modulars.

This research can be considered as a initial step in a broader investigation about quasi-pseudometric modular aggregation functions. While our focus in this work has been primarily on the theoretical aspects, these results pave the way for future applications. In particular, we plan to explore potential applications of this theory in computer science or cluster analysis, where aggregation functions play a critical role. Furthermore, we aim to investigate the relevance of this theory in decision theory, where the use of quasi-pseudometric modulars could provide novel insights into multi-criteria decision-making processes and preference modeling.

Future research will also focus on extending the applicability of the theory of aggregation functions within enriched categories to large classes of mathematical structures and analyzing their practical implications.

\section*{Acknowledgement}{We sincerely thank the anonymous reviewers for their valuable comments and suggestions, which have improved the quality of this paper.}

\end{document}